\documentclass[a4paper,12pt]{article}
\usepackage{amsthm,amsmath}
\usepackage{amssymb}
\usepackage[all]{xypic}
\usepackage{enumerate}
\usepackage{a4wide}
\usepackage{amsthm}
\usepackage{hyperref}

\newtheorem{theorem}{Theorem}[section]
\newtheorem{proposition}[theorem]{Proposition}
\newtheorem{corollary}[theorem]{Corollary}
\newtheorem{lemma}[theorem]{Lemma}
\newtheorem*{theorem*}{Theorem}
\newtheorem*{proposition*}{Proposition}
\newtheorem*{corollary*}{Corollary}
\newtheorem*{lemma*}{Lemma}\theoremstyle{definition}
\newtheorem{definition}[theorem]{Definition}

\newtheorem{example}[theorem]{Example}
\newtheorem{remark}[theorem]{Remark}

\newtheorem*{definition*}{Definition}

\newcommand{\rcomod}[1]{\mathcal{M}^{#1}}

\newcommand{\rmod}[1]{\mathcal{M}_{#1}}

\newcommand{\cotensor}[1]{\square_{#1}}

\renewcommand{\hom}[3]{\mathrm{Hom}_{#1}(#2,#3)}

\newcommand{\cohom}[3]{\mathrm{h}_{#1}(#2,#3)}

\newcommand{\Ext}[4]{\mathrm{Ext}^{#1}_{#2}(#3,#4)}

\newcommand{\cat}[1]{\mathcal{#1}}
\newcommand{\com}[3]{\mathrm{Com}_{#1}(#2,#3)}
\newcommand{\ilimt}[1]{\underset{\longleftarrow}{\mathrm{lim}}\;#1}
\newcommand{\dlimt}[1]{\underset{\longrightarrow}{\mathrm{lim}}\;#1}

\date{}

\begin{document}
\title{Localization in coalgebras. Applications to finiteness conditions.}
\author{J. G\'omez-Torrecillas \footnote{Supported by the
Grant BFM2001-3141 from the Ministerio de Ciencia y Tecnolog\'{\i}a and
FEDER} \\ {\normalsize Departamento de \'{A}lgebra} \\
{\normalsize Facultad de Ciencias}\\ {\normalsize Universidad de Granada} \\
{\normalsize E18071 Granada, Spain} \\ {\normalsize \sf E-mail: torrecil@ugr.es}  \and C. N\u
ast\u asescu \\ {\normalsize University of
Bucharest} \\ {\normalsize Faculty of Mathematics} \\ {\normalsize Str. Academiei 14} \\
{\normalsize RO70109, Bucharest, Romania} \\ {\normalsize \sf E-mail:
cnastase@al.math.unibuc.ro} \and B. Torrecillas \\ {\normalsize Departamento de \'{A}lgebra y
An\'{a}lisis Matem\'{a}tico}
\\ {\normalsize Universidad de Almer\'{\i}a} \\ {\normalsize E04071 Almer\'{\i}a, Spain} \\ {\normalsize \sf E-mail: btorreci@ual.es}}
\maketitle

\section*{Introduction}

Let $k$ be a field. Given two finite-dimensional right comodules $N$ and $M$ over a
$k$--coalgebra $C$, the $k$--vector spaces $Ext_C^n(N,M)$ need not to be finite-dimensional.
This is due to the fact that the injective right comodules appearing in the minimal injective
resolution of $M$ need not to be of finite dimension or even quasi-finite. The obstruction
here is that factor comodules of quasi-finite comodules are not in general quasi-finite. This
note is mainly devoted to the study of coalgebras for which the class of all quasi-finite
right comodules is closed under factor comodules. A major tool here is the local study, in the
sense of abstract localization \cite{Gabriel:1962}, of the comodules which have all their
factors quasi-finite.

\section{Localization in coalgebras}\label{localization}

Let $C$ be a coalgebra over a field $k$, with comultiplication
$\Delta : C \rightarrow C \otimes C$ and counit $\epsilon : C
\rightarrow k$. We will use a variation of Heynemann-Sweedler's
sigma-notation, namely, $\Delta (c) = \sum c_1 \otimes c_2$. The
notation $\rcomod{C}$ stands for the $k$--linear category of all
right $C$--comodules. The hom bifunctor in this category is
denoted by $\com{C}{-}{-}$. Associated to every localizing
subcategory $\cat{T}$ of $\rcomod{C}$ we have the \emph{quotient
category} $\rcomod{C}/\cat{T}$, a $k$--abelian category determined
up to equivalence by an exact functor $T : \rcomod{C} \rightarrow
\rcomod{C}/\cat{T}$ whose kernel is precisely $\cat{T}$, and its
right adjoint $S : \rcomod{C}/\cat{T} \rightarrow \rcomod{C}$; the
adjunction being such that its counit is an isomorphism \cite[Ch.
III]{Gabriel:1962}. By \cite{Takeuchi:1977} there is a
$k$--coalgebra $D$ (unique up to a Morita-Takeuchi equivalence)
such that $\rcomod{C}/\cat{T}$ is a category equivalent to
$\rcomod{D}$. We will give an explicit description of the
coalgebra $D$.

Consider the dual algebra $C^* = \hom{k}{C}{k}$ with the convolution product $f * g (c) = \sum
f(c_1)g(c_2)$. Recall from \cite[Lemma 6]{Radford:1982} that every idempotent $e = e^2 \in
C^*$ gives a $k$--coalgebra $eCe$ with comultiplication defined as $\Delta (ece) = \sum e c_1
e \otimes e c_2 e$. Its counit is the restriction of $\epsilon$. By \cite[Lemma
1.2]{Cuadra/Gomez:2002} we have a $C-eCe$--bicomodule $eC$ (resp. an $eCe-C$--bicomodule
$Ce$), with structure maps defined in a straightforward way. This leads (see \cite[Theorem
1.5, Corollary 1.6]{Cuadra/Gomez:2002}) to an exact functor $- \cotensor{C} eC : \rcomod{C}
\rightarrow \rcomod{eCe}$ with a right adjoint $ - \cotensor{eCe} eC : \rcomod{eCe}
\rightarrow \rcomod{C}$. In fact the functor $- \cotensor{C} eC $ is naturally isomorphic to
the co-hom functor $\cohom{C}{Ce}{-}$ and also to the functor $e(-)$ that sends $M \in
\rcomod{C}$ onto $eM$. The fundamental properties of co-hom functors may be found in
\cite{Takeuchi:1977}. The counit of this adjunction is, by \cite[Proposition
1.4]{Cuadra/Gomez:2002}, an isomorphism, so $\rcomod{eCe}$ becomes a quotient category of
$\rcomod{C}$. By \cite[Proposition 3.8]{Woodcock:1997}, every quotient category of
$\rcomod{C}$ is of this form. We will give here an alternative approach, more appropriate for
our purposes, to this fact.

Let $\{ S_i \; : i \in I \}$ be a complete set of representatives
of the isomorphism types of simple right $C$--comodules, which we
fix from now on. For each subset $J \subseteq I$ we denote by
$\cat{A}_J$ the smallest localizing subcategory of $\rcomod{C}$
that contains the set $\{ S_j \; : \; j \in J \}$. Clearly, $M \in
\cat{A}_J$ if and only if for every subcomodule $M' \lneqq M$, the
quotient $M/M'$ contains a simple subcomodule $S$ with $S \cong
S_j$ for some $j \in J$. Since $\rcomod{C}$ is a locally finite
category, every localizing subcategory $\cat{T}$ of $\rcomod{C}$
is of the form $\cat{A}_J$, where $J \subseteq I$ consists of
those $j \in I$ such that $S_j \in \cat{T}$.

Let $\{ e_i \; : \; i \in I \}$ be a basic set of idempotents for $C$, that is, it is a set of
orthogonal idempotents of $C^*$ such that $E(S_i) = Ce_i$ for every $i \in I$ (see
\cite[Section 3]{Cuadra/Gomez:2002}). The notation $E(M)$ stands for an injective hull (in the
category of comodules) of a comodule $M$. Most part of the following proposition was
essentially given, with a different proof, in \cite[Proposition 3.8]{Woodcock:1997}.

\begin{proposition}\label{locaidem}
Let $\cat{T}$ be a localizing subcategory of $\rcomod{C}$ and let
$J \subseteq I$ such that $\cat{T} = \cat{A}_J$. Let $e = \sum_{i
\notin J}e_i$ be the idempotent acting as zero on $Ce_j$ for $j
\in J$, and as $e_i$ on $Ce_i$ for $i \notin J$.  Then
\begin{enumerate}[(1)]
\item $\cat{T} = \{ M \in \rcomod{C} \; : \; e_iM = 0 \, \hbox{ for all } i \notin J
\} = \{ M \in \rcomod{C} \; : \; eM = 0 \}$.
 \item $\rcomod{C}/\cat{T} = \rcomod{eCe}$ and $- \cotensor{C} eC : \rcomod{C} \rightarrow \rcomod{eCe }$ is the
localizing functor with right adjoint $- \cotensor{eCe} Ce :
\rcomod{eCe} \rightarrow \rcomod{C}$. The localizing functor is
naturally isomorphic to the co-hom functor $\cohom{C}{Ce}{-}$ and
also to the functor $e(-)$ that sends $M \in \rcomod{C}$ onto
$eM$.
\end{enumerate}
\end{proposition}
\begin{proof}
(1) The injective right $C$--comodule $\bigoplus_{i \notin J}
E(S_i)$ is an injective cogenerator for the class of all
$\cat{T}$--torsionfree comodules, that is, $M \in \cat{T}$ if and
only if $\com{C}{M}{\bigoplus_{i \notin J}E(S_i)} = 0$, and this
is equivalent to $\com{C}{M}{Ce_i} = 0$ for every $i \notin J$, as
$E(S_i) = Ce_i$ for every $i$. On the other hand, by \cite[Theorem
1.5]{Cuadra/Gomez:2002}, the co-hom functor $\cohom{C}{Ce_i}{-} :
\rcomod{C} \rightarrow \rcomod{e_iCe_i}$ is naturally isomorphic
to the functor $e_i(-) : \rcomod{C} \rightarrow \rcomod{e_iCe_i}$
for every $i$, where this last functor sends $M \in \rcomod{C}$
onto $e_iM$. Every right $C$--comodule $M$ may be written as a
direct limit $M = \dlimt{M_{\alpha}}$, where the $M_{\alpha}'s$
are finite-dimensional subcomodules. We have
\begin{multline*}
\com{C}{M}{Ce_i} \cong \ilimt{\com{C}{M_{\alpha}}{Ce_i}} \cong
\ilimt{\cohom{C}{Ce_i}{M_{\alpha}}^*} \\ \cong (\dlimt{\cohom{C}{Ce_i}{M_{\alpha}}})^*  \cong
\cohom{C}{Ce_i}{M}^* \cong (e_iM)^*
\end{multline*}
Therefore, $M \in \cat{T}$ if and only if $e_iM = 0$ for all $i
\notin J$ if and only if $eM = 0$. \\
(2) First, note that $eC$ is coflat as a left $C$--comodule, as it
is a direct summand of $C$. Thus, the functor $- \cotensor{C} eC :
\rcomod{C} \rightarrow \rcomod{eCe}$ is exact. By \cite[Corollary
1.6]{Cuadra/Gomez:2002}, $- \cotensor{C} eC$ is naturally
isomorphic to the co-hom functor $\cohom{C}{Ce}{-}$, which means
that $ - \cotensor{eCe}{Ce} : \rcomod{eCe} \rightarrow \rcomod{C}$
is its right adjoint. Now, \cite[Theorem 1.5]{Cuadra/Gomez:2002}
says that $M \cotensor{C} eC$ is naturally isomorphic $eM$ for
every right $C$--comodule. Hence, by part (1), $\cat{T}$ consists
of those comodules $M$ such that $M \cotensor{C} eC = 0$. Finally,
if $N \in \rcomod{eCe}$, then, by \cite[Proposition
1.4]{Cuadra/Gomez:2002}, we have $(N \cotensor{eCe} Ce)
\cotensor{C} eC = N \cotensor{eCe} (Ce \cotensor{C} eC) \cong N
\cotensor{eCe} eCe \cong N$.
\end{proof}

When applied to the case of the localizing subcategory
$\cat{T}_{E(S)}$ consisting of all right comodules $M$ such that
$\com{C}{M}{E(S)} = 0$, where $S$ is a simple right comodule,
Proposition \ref{locaidem} gives a remarkable consequence. Recall
\cite[Definition 1.19]{Cuadra/Gomez:2002} that a coalgebra is said
to be \emph{colocal} if its coradical is the dual of a division
algebra; equivalently, its dual algebra is a local algebra. Every
colocal coalgebra is irreducible, i. e., there it has unique type
of simple comodule.

\begin{corollary}\label{locatsimple}
Let $S$ be a simple right $C$--comodule. If $i \in I$ is such that $S \cong S_i$, then
$\cat{T}_{E(S)} = \cat{A}_{I\setminus \{i\}}$ and  $\rcomod{C}/\cat{T}_{E(S)}$ is equivalent
to $\rcomod{e_iCe_i}$. Moreover, $e_iCe_i$ is a colocal coalgebra and the localizing functor
is $T_i = - \cotensor{C} e_iC : \rcomod{C} \rightarrow \rcomod{e_iCe_i}$.
\end{corollary}

Simple comodules $S$ having a projective cover lead to the
following alternative description of the localization at
$\cat{T}_{E(S)}$.

\begin{proposition}\label{conproy}
Let $S$ be a simple right $C$--comodule and assume that $S$ has a projective cover $P$ (e. g.
$C$ is right semiperfect). Then
\[
\cat{T}_{E(S)} = \{ M \in \rcomod{C} \; | \; \com{C}{P}{M}  = 0 \}
\]
and $\rcomod{C}/\cat{T}_{E(S)}$ is equivalent to $\rmod{R}$, where $R = \com{C}{P}{P}$ is the
endomorphism ring of $P$ in $\rcomod{C}$. Hence, $R$ is a finite-dimensional local ring, and
the colocal coalgebra $e_iCe_i$ given in Corollary \ref{locatsimple} is isomorphic to $R^*$.
\end{proposition}
\begin{proof}
Let $i \in I$ such that $S \cong S_i$. By Corollary
\ref{locatsimple}, $\cat{T}_{E(S)} = \cat{A}_{I \setminus \{i\}}$,
the smallest localizing category containing all $S_j$ with $j \neq
i$. Let $j \in I\setminus \{ i \}$ and assume $\com{C}{P}{S_j}
\neq 0$. Then there is an epimorphism $g : P \rightarrow S_j$.
Since $P$ is local, this would imply that there exists an
epimorphism from $S_j$ onto $S$. Therefore, $S_j \cong S$, a
contradiction. Hence, $\com{C}{P}{S_j} = 0$ for every $j \neq i$,
which implies that $\cat{A}_{I \setminus \{ i \} } \subseteq \{ M
\in \rcomod{C} \; | \; \com{C}{P}{M} = 0\}$. Conversely, let $M$
be a right $C$--comodule such that $\com{C}{P}{M} = 0$, and let
$M'$ be the largest subcomodule of $M$ such that $M' \in
\cat{A}_{I \setminus \{ i \}}$. If $M/M' \neq 0$, then it contains
a simple subcomodule isomorphic to $S_i$ and, hence, to $S$. So,
$\com{C}{P}{M/M'} \neq 0$ and, since $P$ is projective, this
implies that $\com{C}{P}{M} \neq 0$, a contradiction. Therefore,
$M/M' = 0$ and $M = M' \in \cat{A}_{I \setminus \{ i \}}$. Arguing
as in \cite[Proposition 8.6]{Albu/Nastasescu:1984}, we have that
$T(P)$ is a projective generator of the quotient category
$\rcomod{C}/\cat{T}_{E(S)}$, where $T$ denotes the localization
functor. Moreover, $P$ is finite-dimensional, whence $T(P)$ is of
finite length. Therefore, $\rcomod{C}/\cat{T}_{E(S)}$ is
equivalent to the category of left modules over the endomorphism
ring $\hom{\rcomod{C}/\cat{T}_{E(S)}}{T(P)}{T(P)} \cong
\com{C}{P}{P}$.
\end{proof}

\section{Strictly quasi-finite comodules}\label{sqf}

A right comodule $M$ over a coalgebra $C$ is said to be
\emph{quasi-finite} \cite{Takeuchi:1977} if $\com{C}{S}{M}$ is a
finite-dimensional $k$--vector space for every simple right
$C$--comodule $S$ or, equivalently, $\com{C}{N}{M}$ is
finite-dimensional for every comodule $N$ of finite dimension.
Every subcomodule and every essential extension of a quasi-finite
comodule is quasi-finite. However, factor comodules of
quasi-finite comodules are not in general quasi-finite.

\begin{example}\label{ejemplo1}
Let $V$ be a vector space over $k$, and consider $C = kg \oplus V$
be the co-commutative $k$--coalgebra with comultiplication given
by $\Delta(g) = g \otimes g$ and $\Delta(x) = g \otimes x + x
\otimes g$ for $x \in V$, and counit defined by $\epsilon(g) = 1$
and $\epsilon (x) = 0$ for $x \in V$. Its coradical is given by
$C_0 = kg$ and $C = C_0 \wedge C_0$. Hence, $C$ has a unique
simple comodule $kg$ and $C$ is a colocal coalgebra. Moreover,
$C/C_0 \cong V$ is a semisimple comodule. Therefore, if $V$ is not
finite-dimensional, then $C/C_0$ is not a quasi-finite right
$C$--comodule.
\end{example}

\begin{definition}
A right comodule $M$ is said to be \emph{strictly quasi-finite} if
$M/M'$ is quasi-finite for every subcomodule $M'$ of $M$. The
comodule $M$ is said to be \emph{co-noetherian} if $M/M'$ embeds
in a finite direct sum of copies of $C_C$ for every subcomodule
$M'$.
\end{definition}

Co-noetherian comodules were investigated in \cite{Wang/Wu:1998}
and \cite{Wang:1998}. In \cite[Definition 3.2]{Wang:1998},
strictly quasi-finite comodules are also called co-notetherian.

The coalgebra itself is a quasi-finite right comodule. As a
consequence, every co-noetherian comodule is strictly
quasi-finite. The converse is not true in general, as the
following example shows.

\begin{example}\label{ejemplo2}
Let $(V_n)_{n \geq 1}$ be a sequence of $k$--vector spaces such that $dim_k(V_n) = d_n$ with
$d_1 < d_2 < \cdots$. Consider the coalgebras $C_n = kg_n \oplus V_n$ as in Example
\ref{ejemplo1}. Since $C_n$ is finite dimensional for every $n \geq 1$, we have in particular
that $C = \oplus_{n \geq 1} C_n$ is a semiperfect coalgebra. We will see (Theorem
\ref{sqfsemiperfect}) that $C$ is then strictly quasi-finite as a $C$--comodule. Let us show
that $C$ is not co-noetherian as a comodule. Indeed, if $C_C$ were co-noetherian, we had an
exact sequence $0 \rightarrow C/C_0 \rightarrow C^t$ for some $t \geq 1$. From this, we have
an exact sequence of $C_n$--comodules $0 \rightarrow C_n/(C_n)_0 \rightarrow C_n^t$ for every
$n \geq 1$. Hence, $dim_k(\com{C_n}{kg_n}{C_n/(C_n)_0}) \leq dim_k(\com{C_n}{kg_n}{C_n^t}) =
t$ for every $n \geq 1$. Since $C_n/(C_n)_0 \cong V_n$ is semi-simple as a $C_n$--comodule, we
get that $d_n \leq t$ for every $n \geq 1$, which is not possible by the choice of the
dimensions $d_n$. Note that this example also shows that an infinite direct sum of
co-noetherian coalgebras need not to be co-noetherian.
\end{example}

We will say that a class of comodules $\cat{C}$ is a \emph{Serre
class} whenever for any short exact sequence $0 \rightarrow M'
\rightarrow M \rightarrow M'' \rightarrow 0$ of comodules one has
that $M \in \cat{C}$ if and only if $M', M'' \in \cat{C}$. The
first statement of the following proposition was given without
proof in \cite[Proposition 3.1]{Wang:1998}. The second and third
statements were observed in \cite{Cuadra:2001}. We give proofs for
the convenience of the reader. We recall from
\cite{Heyneman/Radford:1974} that a coalgebra $C$ is said to be
\emph{almost connected} if its coradical $C_0$ is
finite-dimensional.

\begin{proposition}\label{sqfyconoeth}\cite{Wang:1998,Cuadra:2001}.
Let $C$ be a coalgebra.
\begin{enumerate}[(1)]
\item The class of all strictly quasi-finite right $C$--comodules
is a Serre class in $\rcomod{C}$.
\item The class of all
co-noetherian right $C$--comodules is a Serre subclass of the
class of all strictly quasi-finite right $C$--comodules.
\item If
$C$ is almost connected, then both classes coincide.
\end{enumerate}
\end{proposition}
\begin{proof}
(1) Let
\begin{equation}\label{ecua1}
\xymatrix{0 \ar[r] & M' \ar[r] & M \ar[r] & M'' \ar[r] & 0}
\end{equation}
be an exact sequence in $\rcomod{C}$. Clearly, if $M$ is strictly
quasi-finite, then $M'$ and $M''$ are strictly quasi-finite.
Assume now that $M'$ and $M''$ are strictly quasi-finite, and let
$X \leq M$ be any subcomodule. Consider the diagram with exact
rows
\begin{equation}\label{ecua2}
\xymatrix{ & & & M'' = M/M' \ar[d]  & \\
           0 \ar[r] & M'/X \cap M' \ar[r] & M/X \ar[r] & M/(M'+X) \ar[r] \ar[d] & 0 \\
           & & & 0 &  }
\end{equation}
By hypothesis, $M'/X \cap M'$ and $M / (M'+X)$ are quasi-finite.
Therefore, we can assume in \eqref{ecua2} that $X = 0$ and we have
to prove that $M$ is quasi-finite in \eqref{ecua1}. This follows
easily from the definition of quasi-finite comodule. \\
(2) Assume in the sequence \eqref{ecua1} that $M$ is
co-noetherian. Clearly, $M'$ is co-noetherian. Now, if $X \leq M'$
is any subcomodule, then there exists $t \geq 0$ such that $M/X$
embeds in $C^t$. Since $M'/X$ is a subcomodule of $M/X$ we get
that $M'/X$ embeds in $C^t$ and so $M'$ is co-noetherian. Assume
now that $M'$ and $M''$ are co-noetherian and let $X \leq M$ be
any subcomodule. From \eqref{ecua2} we can assume that $X = 0$. We
have exact sequences $0 \rightarrow M' \rightarrow C^r$ and $0
\rightarrow M'' \rightarrow C^s$ for some $r, s \geq 0$. Using
that $C^r$ is an injective comodule, we can construct an exact
diagram
\[
\xymatrix{& 0 \ar[d] & 0 \ar[d] & 0 \ar[d] & \\
           0 \ar[r] & M' \ar[r] \ar[d]& M \ar[r] \ar[d] & M'' \ar[r] \ar[d] & 0 \\
          0 \ar[r] & C^r \ar[r] & C^r \oplus C^s \ar[r] & C^{s} \ar[r] &  0}
\]
(3) Over an almost connected coalgebra $C$, every quasi-finite
comodule embeds in a finite direct sum of copies of $C$.
\end{proof}

\begin{proposition}\label{sqfvconoeth}
Let $C$ be an almost connected coalgebra. The following statements
are equivalent for a right $C$--comodule $M$.
\begin{enumerate}[(i)]
\item $M$ is strictly quasi-finite, \item $M$ is co-noetherian,
\item $M^*$ is a noetherian right $C^*$--module, \item $M$ is an
artinian object of $\rcomod{C}$.
\end{enumerate}
\end{proposition}
\begin{proof}
$(i) \Leftrightarrow (ii)$ follows from Proposition
\ref{sqfyconoeth}.\\
$(i) \Rightarrow (iv)$ 
Let $X_1 \supseteq X_2 \supseteq \cdots \supseteq X_n \supseteq \cdots$ be a descending chain
of subcomodules of $M$, and put $X = \bigcap_{i = 1}^{\infty} X_i$. Since $M/X$ is
quasi-finite and $C$ is almost connected, it follows that the socle of $M/X$, which is
essential, is finite-dimensional. Arguing as in \cite[Proposition
10.10]{Anderson/Fuller:1992}, there is a positive integer $n$ such that
$X = \bigcap_{i = 1}^n X_i$. Hence, $X_n = X_{n+1} = \cdots $. \\
$(iii) \Rightarrow (iv)$ With notations as above, if $X_i^{\bot} =
\{ f \in M^* \; | \; f(X_i) = 0 \}$ then we have an ascending
chain $X_1^{\bot} \subseteq X_2^{\bot} \subseteq \cdots \subseteq
X_n^{\bot} \subseteq \cdots$ of $C^*$--submodules of $M^*$ . Since
$M^*$ is noetherian, we get $X_n^{\bot} = X_{n+1}^{\bot} = \cdots
$ for some positive integer $n$. But $(X_n/X_{n+1})^* \cong
X_{n+1}^{\bot}/X_n^{\bot} = 0$, whence
$X_n/X_{n+1} = 0$. \\
 $(iv) \Leftrightarrow (iii)$ follows from
\cite[Corollary 4.3]{Dascalescu/Nastasescu/Raianu:2001}.\\
$(iv) \Rightarrow (i)$ Since $M$ is artinian, every factor
comodule $M/X$ of $M$ has finite-dimensional socle and, in
particular, $M/X$ is quasi-finite.
\end{proof}

We will see later (Example \ref{noartiniano}) that a strictly quasi-finite comodule over an
arbitrary coalgebra need not to be artinian. It shares, however, some properties with artinian
objects. An example is the following.

\begin{proposition}
Let $M$ be a strictly quasi-finite right comodule over a coalgebra
$C$. If $f : M \rightarrow M$ is a monomorphism in $\rcomod{C}$
then $f$ is an isomorphism.
\end{proposition}
\begin{proof}
Let $M_0 \subset M_1 \subset \cdots $ be the Loewy series of $M$,
that is, $M_{n+1}/M_n$ is the socle of $M/M_n$ for $n > 0$ and
$M_0$ is the socle of $M$. It suffices to prove that $f(M_n) =
M_n$ for every $n \geq 0$, as $M = \bigcup_{n \geq 0}M_n$. On the
other hand, $I = C_0^\bot$ is the Jacobson radical of $C^*$, and
$M_n = ann_M(I^{n+1})$ for every $n \geq 0$ (see \cite[Lemma
3.1.9]{Dascalescu/Nastasescu/Raianu:2001}). We proceed by
induction on $n$. For $n = 0$ we have that $f(M_0) \subseteq M_0$.
Since $M$ is quasi-finite, every isotypic component of $M$ is
finite-dimensional, which gives $f(M_0) = M_0$. Assume inductively
that $f(M_n) = M_n$. Then we have $f(M_{n+1}) \subseteq M_{n+1}$,
and we can consider the induced morphism $\overline{f} :
M_{n+1}/M_n \rightarrow M_{n+1}/M_n$. If $\overline{f}( x + M_n) =
0$, then $f(x) \in M_n$, so $I^{n+1}f(x) = 0$. Hence, $f(I^{n+1}x)
= 0$ which implies that $I^{n+1}x = 0$, as $f$ is injective.
Therefore, $x \in M_n$ and $\overline{f}$ is injective. This
implies that $\overline{f}$ is bijective because $M/M_n$ is
quasi-finite. So, $f(M_{n+1}) = M_{n+1}$ which completes the
induction.
\end{proof}

Our next aim is to characterize strictly quasi-finite comodules in
terms of their localizations at the simple comodules.

\begin{lemma}\label{locasimple}
Let $\cat{T}$ be any localizing subcategory of the category
$\rcomod{C}$ of right comodules over a coalgebra $C$, and let $T :
\rcomod{C} \rightarrow \rcomod{C}/\cat{T}$ be the localizing
functor. If $X$ is a simple object in the $\rcomod{C}/\cat{T}$,
then there exists a simple right $C$--comodule $Y$ such that $T(Y)
\cong X$.
\end{lemma}
\begin{proof}
Let $S : \rcomod{C}/\cat{T} \rightarrow \rcomod{C}$ be the right
adjoint functor to $T$. Then $S(X)$ is a nonzero right
$C$--comodule, so there is an exact sequence $0 \rightarrow Y
\rightarrow S(X)$, where $Y$ is a simple comodule. Since $T$ is
exact, we get the exact sequence $0 \rightarrow S(Y) \rightarrow
TS(X) \cong X$, which gives an isomorphism $T(Y) \cong X$, as $X$
is a simple object.
\end{proof}

\begin{proposition}\label{MsqfTMartinian}
Let $\cat{T}$ be a localizing subcategory of the category
$\rcomod{C}$ of right comodules over a coalgebra $C$, and let $T :
\rcomod{C} \rightarrow \rcomod{C}/\cat{T}$ be the localization
functor.
\begin{enumerate}[(1)]
\item If $M$ is a strictly quasi-finite right $C$--comodule, then
$T(M)$ is strictly quasi-finite. \item If $\cat{T} =
\cat{T}_{E(S)}$ for some simple right $C$--comodule $S$ and $M$ is
a strictly quasi-finite right $C$--comodule, then $T(M)$ is an
artinian object of $\rcomod{C}/\cat{T}_{E(S)}$.
\end{enumerate}
\end{proposition}
\begin{proof}
(1) Since $T$ is an exact functor, it suffices to prove that
$T(M)$ is quasi-finite. Let $X$ be a simple object of
$\rcomod{C}/\cat{T}$. By Lemma \ref{locasimple}, there is a simple
right $C$--comodule $Y$ such that $T(Y) \cong X$. Thus,
\[
\hom{\rcomod{C}/\cat{T}}{X}{T(M)} \cong
\hom{\rcomod{C}/\cat{T}}{T(Y)}{T(M)} \cong \com{C}{Y}{ST(M)} ,
\]
where $S : \rcomod{C}/\cat{T} \rightarrow \rcomod{C}$ denotes the
right adjoint to $T$. The kernel of the canonical map $\alpha : M
\rightarrow ST(M)$ is an object of $\cat{T}$, which gives an exact
sequence of vector spaces $0 \rightarrow \com{C}{Y}{M} \rightarrow
\com{C}{Y}{Im (\alpha})$. On the other hand, $Im (\alpha)$ is
essential in $ST(M)$ and, since $Im (\alpha)$ is quasi-finite, we
have that $ST(M)$ is quasi-finite too. Therefore,
$\com{C}{Y}{ST(M)}$ is finite-dimensional and so is
$\hom{\rcomod{C}/\cat{T}}{X}{T(M)}$. \\
(2) This follows from Corollary \ref{locatsimple} and Proposition
\ref{sqfvconoeth} in conjunction with part (1).
\end{proof}

Let $\{ S_i \; : i \in I \}$ be a complete set of representatives of the isomorphism types of
simple right $C$--comodules. For every $i \in I$ let $T_i : \rcomod{C} \rightarrow
\rcomod{C}/\cat{T}_{E(S_i)}$ denote the localization functor.

\begin{theorem}\label{MsqfiffMiartinian}
A right $C$--comodule $M$ is strictly quasi-finite if and only if
$T_i(M)$ is artinian for every $i \in I$.
\end{theorem}
\begin{proof}
If $M$ is quasi-finite, then $T_i(M)$ is artinian for every $i \in
I$ by Proposition \ref{MsqfTMartinian}. Conversely, assume
$T_i(M)$ artinian for every $i \in I$. Since the functors $T_i$
are exact, it is enough to prove that $M$ is quasi-finite. Let $S$
be a simple right $C$--comodule and consider the unique $i \in I$
such that $S \cong S_i$. If $M_i$ denotes the largest subcomodule
of $M$ such that $M_i \in \cat{T}_{E(S_i)}$, then we have an exact
sequence $0 \rightarrow M_i \rightarrow M \rightarrow M/M_i
\rightarrow 0$. Then $\com{C}{S_i}{M_i} = 0$ and $\com{C}{S_i}{M}$
is isomorphic to a vector subspace of $\com{C}{S_i}{M/M_i}$. Thus,
we can assume that $M_i = 0$, that is, $M$ is
$\cat{T}_{E(S_i)}$--torsionfree. The socle $soc(M)$ is then a
direct sum of copies of $S_i$ and, since $T_i(M)$ is artinian, we
have that $T_i(soc(M)) \subseteq T_i(M)$ is artinian as well.
Therefore, $soc(M)$ consists of a direct sum of finitely many
copies of $S_i \cong S$ and, hence, $\com{C}{S}{M}$ is
finite-dimensional.
\end{proof}

Recall from \cite{Lin:1977} that a coalgebra is said to be \emph{right semiperfect} if every
simple right comodule has a projective cover.

\begin{theorem}\label{sqfsemiperfect}
Let $C$ be a right semiperfect coalgebra. The following statements are equivalent for a right
$C$--comodule $M$.
\begin{enumerate}[(i)]
\item $M$ is quasi-finite; \item $M$ is strictly quasi-finite;
\item $T_i(M)$ is finite-dimensional for every $i \in I$.
\end{enumerate}
\end{theorem}
\begin{proof}
$(i) \Rightarrow (ii)$ Let $S$ be a simple right $C$--comodule. The projective cover $P
\rightarrow S \rightarrow 0$ of $S$ gives, for every subcomodule $M' \leq M$, two exact
sequences of vector spaces $0 \rightarrow \com{C}{S}{M/M'} \rightarrow \com{C}{P}{M/M'}$ and
$\com{C}{P}{M} \rightarrow \com{C}{P}{M/M'} \rightarrow 0$. Therefore, $\com{C}{S}{M/M'}$ is
finite-dimensional, as $P$ has finite dimension and $M$ is quasi-finite. This proves that
$M/M'$ is quasi-finite and, so, $M$ is strictly quasi-finite. \\
$(ii) \Rightarrow (iii)$ This follows from Proposition \ref{conproy} and Theorem
\ref{MsqfiffMiartinian}.\\
$(iii) \Rightarrow (i)$ By Theorem \ref{MsqfiffMiartinian}. 
\end{proof}

We can now give an alternative proof of one of the most useful characterizations of right
semiperfect coalgebras from \cite{Lin:1977}.

\begin{corollary}\label{Lin}\textbf{(Lin)}
A coalgebra $C$ is right semiperfect if and only if $T_i(C) (= e_iC)$ is finite-dimensional
for every $i \in I$.
\end{corollary}
\begin{proof}
If $C$ is right semiperfect then, by Theorem \ref{sqfsemiperfect}, $T_i(C) = e_iC$ is
finite-dimensional for every $i \in I$, as $C$ is quasi-finite. Conversely, given any simple
right $C$--comodule $S$, pick $i \in I$ such that $E(S^*) = e_iC$. Then the essential
inclusion $S^* \subseteq e_iC$ gives an epimorphism with small kernel $(e_iC)^* \rightarrow
S^{**} \cong S$, where $(e_iC)^*$ is a right $C$--comodule because $e_iC$ is a
finite-dimensional left $C$--comodule. By \cite[Lemma 1.2]{Cuadra/Gomez:2002}, $(e_iC)^* \cong
C^*e_i$, and therefore it is a projective cover of $S$.
\end{proof}

\section{Strictly quasi-finite coalgebras}

Recall that $\{ S_i \; : \; i \in I \}$  is a complete set of representatives of the
isomorphism types of simple right comodules over a coalgebra $C$. For each $i \in I$, let $T_i
: \rcomod{C} \rightarrow \rcomod{C}/\cat{T}_{E(S_i)}$ denote the canonical localization
functor. The following is the main theorem of this section.

\begin{theorem}
The following statements are equivalent for a coalgebra $C$.
\begin{enumerate}[(i)]
\item $C$ is strictly quasi-finite as a right $C$--comodule, \item
for every $i \in I$, $E(S_i)$ is strictly quasi-finite and
$T_i(E(S_j)) \neq 0$ only for finitely many $j \in I$, \item every
quasi-finite right $C$--comodule is strictly quasi-finite, \item
$T_i(C)$ is an artinian object for every $i \in I$.
\end{enumerate}
\end{theorem}
\begin{proof}
$(i) \Rightarrow (ii)$ Given $i \in I$, $E(S_i)$ embeds in $C$, so
$E(S_i)$ is strictly quasi-finite. Now, decompose the injective
right comodule $C$ as $C = \bigoplus_{j \in I}^{m_j}E(S_j)^{m_j}$,
where the $m_j$ are positive integers. By Theorem
\ref{MsqfiffMiartinian}, $T_i(C) = \bigoplus_{j \in
I}T_i(E(S_j))^{m_j}$ is an
artinan object. Therefore, only finitely many nonzero direct summands should appear. \\
$(ii) \Rightarrow (iii)$ Let $M$ be any quasi-finite right
$C$--comodule and consider $i \in I$. Then $E(M)$ is quasi-finite
and, thus, $E(M) \cong \bigoplus_{j \in I}E(S_j)^{n_j}$ for some
set of non-negative integers $\{ n_j  \; ; \; j \in I \}$.
Therefore, $T_i(E(M)) = \bigoplus_{j \in I}T_i(E(S_j))^{m_j}$ has
finitely many nonzero artinian summands (Theorem
\ref{MsqfiffMiartinian}) and, thus, $T_i(E(M))$ is artinian. Using
Theorem \ref{MsqfiffMiartinian} once more, we obtain that $E(M)$
is
strictly quasi-finite and so is $M$.\\
$(iii) \Rightarrow (i)$ This is clear, as $C_C$ is quasi-finite.\\
$(i) \Leftrightarrow (iv)$ By Theorem \ref{MsqfiffMiartinian}.
\end{proof}

\begin{definition}
A coalgebra $C$ satisfying the equivalent conditions of Theorem
\ref{MsqfiffMiartinian} will be said to be \emph{right strictly
quasi-finite}. Every right semiperfect coalgebra is right strictly
quasi-finite (Theorem \ref{sqfsemiperfect}).
\end{definition}

\begin{remark}
The implication $(i) \Rightarrow (iii)$ was left open in \cite[p.
460]{Wang:1998}.
\end{remark}

Let $\{ e_i \; : \; i \in I \}$ be a complete set of orthogonal primitive idempotents for $C$
(see Section \ref{localization}) . Then right strictly quasi-finite coalgebras can be
characterized as follows.

\begin{corollary}\label{eiej}
The coalgebra $C$ is right strictly quasi-finite if and only if
for every fixed $i \in I$ the right $e_jCe_j$--comodule $e_jCe_i$
is artinian for every $j \in I$ and $e_iCe_j \neq 0$ only for
finitely many $j \in I$.
\end{corollary}

Of course, not every coalgebra is strictly quasi-finite (see
Example \ref{ejemplo1}). The following example shows that the
class of strictly quasi-finite coalgebras contains strictly the
class of semiperfect coalgebras.

\begin{example}\label{kX}
Let $C = k[X]$ the Hopf algebra of polynomials in one
indeterminate $X$; its structure of coalgebra is given by
$\Delta(X) = X \otimes 1 + 1 \otimes X$ and $\epsilon (X) = 0$.
Since $C^* \cong k[[X]]$ is a noetherian algebra, it follows from
\ref{sqfvconoeth} that $C$ is an artinian comodule. Having just
one type of simple, this coalgebra is then strictly quasi-finite.
However, it is not semiperfect.
\end{example}

The notion of a right strictly quasi-finite coalgebra is not
left-right symmetric, as the following example shows.

\begin{example}\cite[Example
3.2.9]{Dascalescu/Nastasescu/Raianu:2001} Let $C$ be the
$k$--coalgebra with basis $\{ g_n, d_n \; : \; n \geq 1 \}$ and
structure maps given by $\Delta (g_n) = g_n \otimes g_n$, $\Delta
(d_n) = g_1 \otimes d_n + d_n \otimes g_{n+1}$, $\epsilon (g_n) =
1$ and $\epsilon (d_n) = 0$. Define $g_n^* \in C^*$ by $g_n^*(g_m)
= \delta_{nm}$ and $g_n^*(s_m) = 0$ for every $m \geq 1$. We have
that $C = \oplus_{n \geq 1} g_n^*C$, and $g_n^*C$ is the injective
envelope of the simple left $C$--comodule $kg_n$ for every $n \geq
1$. Thus, $\{ g_n^* \; : \; n \geq 1 \}$ is a complete set of
primitive ortogonal idempotents for the coalgebra $C$. An easy
computation gives that $g_n^*Cg_1 = kd_{n+1} \neq 0$ for every $n
> 1$ which implies, by Corollary \ref{eiej} that $C$ is not left
strictly quasi-finite. On the other hand, $C$ is right
semiperfect, as $g_n^*C = kg_n \oplus kd_{n-1}$ for $n > 1$ and
$g_1^*C = kg_1$.
\end{example}

The following proposition collects a number of properties of right strictly quasi-finite
coalgebras. Recall that a $k$--coalgebra $D$ is said to be \emph{Morita-Takeuchi equivalent}
to a $k$--coalgebra $C$ if there is a $k$--linear equivalence of categories $F : \rcomod{D}
\rightarrow \rcomod{C}$.

\begin{proposition}\label{propiedades}
\begin{enumerate}[(1)]
\item Every subcoalgebra of a right strictly quasi-finite
coalgebra is right strictly quasi-finite.
\item If $C$ is a right
strictly quasi-finite coalgebra and $A \subseteq C$ is a
finite-dimensional subcoalgebra, then $A_{\infty} = \bigcup_{n
\geq 1} \bigwedge^n A$ is right artinian as a right comodule. When
$C$ is right semiperfect, $A_{\infty}$ becomes finite-dimensional.
\item If $C$ is a right strictly quasi-finite coalgebra and $D$
is a coalgebra Morita-Takeuchi equivalent to $C$, then $D$ is
right strictly quasi-finite.
\item Any direct sum of right
strictly quasi-finite coalgebras is right strictly quasi-finite.
\end{enumerate}
\end{proposition}
\begin{proof}
(1) Easy.\\
(2) Let $\rho_M : M \rightarrow M \otimes C$ denote the structure map of a right $C$--comodule
$M$. If $\cat{C}_A = \{ M \in \rcomod{C} \; ; \; \rho_M(M) \subseteq M \otimes A \}$ denotes
the closed subcategory associated to $A$, then $\cat{C}_{A_{\infty}}$ is the smallest
localizing subcategory of $\rcomod{C}$ containing $\cat{C}_A$ (see
\cite{Nastasescu/Torrecillas:1994}). Clearly, $\cat{C}_{A_{\infty}}$ contains only finitely
many isomorphism types of simple comodules. Thus $A_{\infty}$ is an almost connected
coalgebra, so, by Proposition \ref{sqfvconoeth}, $A_{\infty}$ is artinian as a right comodule.
The statement for the semiperfect case follows from
\cite[Corollary 3.2.11, Exercise 3.3.13]{Dascalescu/Nastasescu/Raianu:2001}.\\
(3) Let $C$ and $D$ be Morita-Takeuchi equivalent coalgebras and
assume $C$ to be right strictly quasi-finite. Let $F : \rcomod{D}
\rightarrow \rcomod{C}$ denote an equivalence of categories. Let
$M \leq D$ be any right subcomodule. Now, $F(D)$ is a quasi-finite
right $C$--comodule, which implies, by Theorem
\ref{MsqfiffMiartinian}, that $F(D/M) = F(D)/F(M)$ is
quasi-finite, as $C$ is right strictly quasi-finite. Therefore,
$D/M$ is quasi-finite, i.e., $D$ is right strictly quasi-finite.
\\
(4) Assume $C = \bigoplus_{\alpha \in \Lambda} C_{\alpha}$, as a coalgebra. By \cite[Exercise
2.2.18]{Dascalescu/Nastasescu/Raianu:2001}, $\rcomod{C} = \prod_{\alpha \in \Lambda}
\rcomod{C_{\alpha}}$, which means that every right $C$--comodule $M$ decomposes uniquely as $M
= \bigoplus_{\alpha \in \Lambda} M_{\alpha}$, where $M_{\alpha} \in \rcomod{C_{\alpha}}$ for
every $\alpha \in \Lambda$. This implies that every factor comodule of $M$ is of the form $M/N
= \bigoplus_{\alpha \in \Lambda} M_{\alpha}/N_{\alpha}$ for appropriate
$C_{\alpha}$--subcomodules $N_{\alpha}$ of the $M_{\alpha}$'s. If $X$ is a simple right
$C$--comodule, then $X \in \rcomod{C_{\alpha}}$ for an uniquely determined $\alpha$. It
follows that $\com{C}{X}{M/N} = \com{C_{\alpha}}{X}{M_{\alpha}/N_{\alpha}}$. The statement
follows now taking $M = C$.
\end{proof}

\begin{example}\label{noartiniano}
Let $C = \bigoplus_{\alpha \in \Lambda}C_{\alpha}$ any infinite
direct sum of right strictly quasi-finite coalgebras. By
Proposition \ref{propiedades}, $C_C$ is a strictly quasi-finite
right comodule, and it is not artinian.
\end{example}

The following consequence of Corollary \ref{eiej} shows that co-commutative strictly
quasi-finite coalgebras are just the direct sums of artinian co-commutative coalgebras.

\begin{corollary}
The following statements are equivalent for a co-commutative
coalgebra $C$.
\begin{enumerate}[(i)]
\item $C$ is strictly quasi-finite;
\item $e_iCe_i$ is artinian
for every $i \in I$;
\item $C$ is a direct sum of (co-commutative) artinian colocal coalgebras;
\item $(e_iCe_i)^*$ is noetherian for every $i \in I$.
\end{enumerate}
\end{corollary}
\begin{proof}
$(i) \Rightarrow (ii)$ This follows from Corollary \ref{eiej}.\\
$(ii) \Rightarrow (iii)$ Since $C$ is co-commutative, we get $C = \oplus_{i \in I} e_iCe_i$.\\
$(iii) \Rightarrow (i)$ Apply Proposition \ref{propiedades}.4.\\
$(ii) \Leftrightarrow (iv)$ By Proposition \ref{sqfvconoeth}.
\end{proof}

\begin{remark}
The equivalence between $(i)$ and $(iv)$ was obtained in \cite{Cuadra:2001} by different
methods.
\end{remark}

\section{Some homological coalgebra}

If $M$ is a quasi-finite right comodule over a right strictly
quasi-finite coalgebra $C$, then all terms in the minimal
injective resolution
\begin{equation}\label{ecua3}
\xymatrix{0 \ar[r] & M \ar[r]^{\varepsilon} & Q_0 \ar[r] & Q_1
\ar[r] & \cdots \ar[r] & Q_n \ar[r] & \cdots}
\end{equation}
are quasi-finite. This is easily deduced from the fact that every
factor of a quasi-finite comodule is quasi-finite (Theorem
\ref{MsqfiffMiartinian}) and that the injective envelope of a
quasi-finite comodule is quasi-finite. Now, every quasi-finite
injective right $C$--comodule $Q$ decomposes as $Q = \bigoplus_{i
\in I}E(S_i)^{\mathsf{n}_{Q,i}}$, where $\mathsf{n}_{Q,i}$ is a
non negative integer for every $i \in I$. By
Azumaya-Krull-Remak-Schmidt's Theorem, the numbers
$\mathsf{n}_{Q,i}$ are uniquely determined by $Q$, and hence they
are invariants of the comodule $Q$. As a consequence, the numbers
$\mathsf{n}_{Q_n,i}$ for $i \in I$ and $n \geq 0$ are invariants
of the quasi-finite comodule $M$. If $S \cong S_i$, then we use
the notation $\mathsf{n}_{Q_n,S}$ to refer to $\mathsf{n}_{Q_n,i}$
They can be computed using $\mathrm{Ext}$--functors, as we shall
show. The notation $\Ext{n}{C}{N}{-}$ stands for the $n$-th
derived functor of $\com{C}{N}{-} : \rcomod{C} \rightarrow
\rmod{k}$, where $N$ is a right $C$--comodule and $n$ is a non
negative integer.

\begin{corollary}\label{Bassnumber}
Let $M$ be a quasi-finite right comodule over a right strictly
quasi-finite coalgebra $C$ with minimal injective resolution as in
\eqref{ecua3}.
\begin{enumerate}[(1)]
\item $\Ext{n}{C}{N}{M}$ is finite-dimensional for every finite-dimensional
right $C$--comodule $N$.
\item If $S$ is a simple right $C$--comodule, then
\[
\mathsf{n}_{Q_n,S} = \frac{dim_k \Ext{n}{C}{S}{M}}{dim_k
\com{C}{S}{S}}
\]
for every $n \geq 0$.
\end{enumerate}
\end{corollary}
\begin{proof}
(1) It is clear, since all terms $Q_n$ are quasi-finite. \\
(2) Consider the short exact sequence
\[
\xymatrix{0 \ar[r] & M \ar[r]^{\varepsilon} & Q_0 \ar[r] & K_1
\ar[r] & 0,}
\]
where $K_1$ is the cokernel of $\varepsilon$. We have an exact
sequence
\[
\xymatrix{0 \ar[r] & \com{C}{S}{M} \ar[r] & \com{C}{S}{Q_0}
\ar[r]& \com{C}{S}{K_1} \ar[r] & \Ext{1}{C}{S}{M} \ar[r] & 0 }
\]
Since $\com{C}{S}{M} = \com{C}{S}{Q_0}$ we get $\com{C}{S}{K_1}
\cong \Ext{1}{C}{S}{M}$. Now, $\com{C}{S}{K_1} = \com{C}{S}{Q_1}$,
as $Q_1 = E(K_1)$, so $\com{C}{S}{Q_1} \cong \Ext{1}{C}{S}{M}$. On
the other hand, $\com{S}{S}{Q_1} =
\com{C}{S}{S}^{\mathsf{n}_{Q_1,S}}$, and therefore $dim_k
\Ext{1}{C}{S}{M} = n_{Q_1,S} \times dim_k \com{C}{S}{S}$. An easy
induction, using the isomorphism $\Ext{n-1}{C}{S}{K_1} \cong
\Ext{n}{C}{S}{M}$ for $n > 1$, gives the desired equality for
every $n$.
\end{proof}

\begin{remark}
If $C$ is right semiperfect, then every projective object of
$\rcomod{C}$ is projective as a left $C^*$--module \cite[Lemma
2.1]{Gomez/Nastasescu:1995}, \cite[Corollary
2.4.22]{Dascalescu/Nastasescu/Raianu:2001}. We thus have
$\Ext{n}{C}{N}{M} = \Ext{n}{C^*}{N}{M}$ for all right
$C$--comodules $N, M$.
\end{remark}

A coalgebra is said to be right \emph{hereditary}
\cite{Nastasescu/Torrecillas/Zhang:1996} if every factor comodule
of an injective right comodule is injective. Our closing result
says that a right hereditary right semiperfect coalgebra is
``locally'' semisimple. An example that shows that the pertinence
of the semiperfect hypothesis is included.

\begin{proposition}
Let $C$ be a right semiperfect right hereditary coalgebra. Then
$e_iCe_i$ is dual to a division $k$--algebra for every $i \in I$.
\end{proposition}
\begin{proof}
Some standard arguments show that if $\cat{A} \rightarrow
\cat{A}/\cat{T}$ is a localization of a Grothendieck category
$\cat{A}$ such that factors of injective objects are injective,
then the quotient category $\cat{A}/\cat{T}$ inherits such a
property. Therefore, if $C$ in a right hereditary coalgebra, then
so is $C_i$ (Proposition \ref{locaidem}). According to Corollary
\ref{Lin}, $C_i = e_iCe_i$ is finite-dimensional for every $i \in
I$. Being a colocal coalgebra \cite[Proposition
1.20]{Cuadra/Gomez:2002}, there is an epimorphism of right
$C$--comodules from $C_i$ onto its coradical $(C_i)_0$, which is
the dual of a (finite dimensional) division algebra. This implies
that  $(C_i)_0$ is an injective right $C_i$--comodule, as $C_i$ is
right hereditary. Therefore, $C_i = (C_i)_0$, which finishes the
proof.
\end{proof}

\begin{example}
The coalgebra $k[X]$ given in Example \ref{kX} is strictly
quasi-finite, colocal and hereditary, by it is obviously not the
dual of a division algebra.
\end{example}

 \providecommand{\bysame}{\leavevmode\hbox to3em{\hrulefill}\thinspace}
\providecommand{\MR}{\relax\ifhmode\unskip\space\fi MR }
\providecommand{\MRhref}[2]{%
  \href{http://www.ams.org/mathscinet-getitem?mr=#1}{#2}
} \providecommand{\href}[2]{#2}

\end{document}